\documentclass[11pt]{article}
\usepackage{amscd}
\usepackage{amsfonts}
\usepackage{amsmath}
\usepackage{amssymb}
\usepackage{amsthm}
\usepackage{bbm}
\usepackage{CJK}
\usepackage{fancyhdr}
\usepackage{graphicx}
\usepackage{hyperref}
\usepackage{indentfirst}
\usepackage{latexsym}
\usepackage{mathrsfs}
\usepackage{xypic}
\usepackage[compress]{cite}

\newtheorem{theorem}{Theorem}[section]
\newtheorem{lemma}[theorem]{Lemma}
\newtheorem{definition}[theorem]{Definition}
\newtheorem{proposition}[theorem]{Proposition}
\newtheorem{example}[theorem]{Example}

\newtheorem{remark}[theorem]{Remark}

\usepackage[top=1in,bottom=1in,left=1.25in,right=1.25in]{geometry}
\def\<{\langle}
\def\>{\rangle}

\def\o{\otimes}

\date{\today}
\begin{document}
\renewcommand{\baselinestretch}{1.2}
\renewcommand{\arraystretch}{1.0}
\title{\bf Cohomologies of  Reynolds Lie-Yamaguti algebras of any weight  and applications}
\date{\today }

\author{{\bf Wen Teng, Shuangjian Guo\footnote
        { Corresponding author:~~shuangjianguo@126.com}  }\\
{\small  School of Mathematics and Statistics, Guizhou University of Finance and Economics} \\
{\small  Guiyang  550025, P. R. of China}}
 \maketitle
\begin{center}
\begin{minipage}{13.cm}

{\bf \begin{center} ABSTRACT \end{center}}
The purpose of the present paper is to investigate cohomologies of Reynolds Lie-Yamaguti algebras of any weight  and  provide some applications. First,  we introduce  the  notion of  Reynolds Lie-Yamaguti algebras and give some new examples. Moreover, cohomologies of  Reynolds operators and Reynolds Lie-Yamaguti algebras with coefficients in a suitable representation are established.  Finally,    formal deformations and  abelian extensions of    Reynolds Lie-Yamaguti algebras   are characterized in terms of lower degree cohomology groups.
 \smallskip

{\bf Key words}: Lie-Yamaguti algebra,  Reynolds operator, representation,  cohomology,
 \smallskip

 {\bf 2020 MSC:}    17B38, 17B60, 17B56,  17D99
 \end{minipage}
 \end{center}
 \normalsize\vskip0.5cm

\section{Introduction}
\def\theequation{\arabic{section}. \arabic{equation}}
\setcounter{equation} {0}

The concept of Rota-Baxter operators on associative algebras was introduced  by  Baxter \cite{B60}
 in his study of fluctuation theory in probability,  and  further studied by Rota \cite{R69} in algebra and combination.
Rota-Baxter operators have been widely studied in  many areas of mathematics and physics, including  combinatorics, number theory, operads and quantum field theory \cite{C00}. For more details on the Rota-Baxter operators, see \cite{G12}.  Later, other operators related to   Rota-Baxter operators   emerged continuously.
One of them is  Reynolds operator,  which is inspired by Reynolds' work on turbulence in fluid dynamics \cite{Reynolds}.
  Kamp$\mathrm{\acute{e}}$ de F$\mathrm{\acute{e}}$riet \cite{Kamp} created the notion of the Reynolds operator   as a mathematical subject in general.
 And then, Zhang,  Gao and Guo gave examples and properties of Reynolds operators, and studied the free Reynolds algebras by using bracketed words
and rooted trees in \cite{ZGG}.  Further research on Reynolds operators has been developed,  see  \cite{Das1,Gharbi,Hou,Wang,T23}.

In recent years, scholars have increasingly focused on operator structures with any weight.  Wang and Zhou \cite{Wang33,Wang22} developed deformations and homotopy theory  for   Rota-Baxter associative
algebras of any weight. Later, Das has considered the cohomology and deformations of weighted Rota-Baxter Lie algebras by use of Wang and Zhou's method  in   \cite{Das34}.  The authors \cite{Guo36,Hou35} established the cohomology, extensions and deformations of Rota-Baxter 3-Lie algebras with any weight.
In  \cite{Chen37}, Chen, Lou and Sun considered the cohomology and  extensions of Rota-Baxter Lie triple systems.

The concept of Lie-Yamaguti algebras  was introduced by Kinyon and Weinstein \cite{K01} in the study of Courant  algebroids.  A Lie-Yamaguti algebra is a generalization of a Lie algebra  and a Lie triple system.
Its roots can be traced back to Yamaguti's study on general Lie triple systems \cite{Y57} and Kikkawa's research on Lie triple algebras   \cite{K75,K79,K81}.
Yamaguti   established the representation and cohomology theory of   Lie-Yamaguti algebras  in  \cite{Y67}.  In \cite{B05,B09},  Benito and his collaborators   studied Lie-Yamaguti algebras related to simple Lie algebras and irreducible Lie-Yamaguti algebras. The  deformation and extension theory for Lie-Yamaguti algebras were investigated in \cite {L15,Z15,G24}.
 The authors of the present paper studied the cohomology theory and deformations of weighted Rota-Baxter Lie-Yamaguti algebras \cite{Teng2023}.
See \cite{Teng2024,TG2024,T2024,S21,S22,Z22,L23} for some interesting related about Lie-Yamaguti algebras.

Motivated by \cite{Wang33,Wang22} and \cite{Das34},   our aim in this paper is to consider   Reynolds Lie-Yamaguti algebras  of any weight, and construct the cohomology theory of   Reynolds Lie-Yamaguti algebras with
coefficients in a suitable representation.
As applications of our cohomology,  we discuss the  formal deformation   and abelian extension  of  Reynolds Lie-Yamaguti algebras.

  This paper is organized as follows.
In Section  \ref{sec:Preliminaries}, we recall some basic concepts about Lie-Yamaguti  algebras,  representations and their cohomology  theory.
In  Section \ref{sec:Representations}, we  consider Reynolds Lie-Yamaguti algebras  of any weight and  present  their representations.
In Section \ref{sec:Cohomology},  we construct the cohomology of a Reynolds Lie-Yamaguti algebra  with coefficients in a suitable representation..
In Section \ref{sec:deformations},  we use the cohomological approach to study formal deformations of Reynolds Lie-Yamaguti algebras.
 Finally, abelian extensions of Reynolds Lie-Yamaguti algebras are studied in Section  \ref{sec:extensions}.

\section{Lie-Yamaguti  algebras, representations  and their cohomology theory }\label{sec:Preliminaries}
\def\theequation{\arabic{section}. \arabic{equation}}
\setcounter{equation} {0}

Throughout this paper, we work on an algebraically closed field $\mathbb{K}$ of characteristic zero.  In this section, we  recall some basic  definitions of  Lie-Yamaguti  algebras from \cite{K01} and \cite{Y67}.

\begin{definition}
A Lie-Yamaguti  algebra  is a  triple  $(L, [-, -], \{-, -, -\})$ in which $L$ is a vector space together with  a binary bracket $[-, -]$ and  a ternary bracket $\{-, -, -\}$
on $L$   satisfying
\begin{eqnarray*}
&&(LY1)~~ [x, y]=- [y,  x],\\
&&(LY2) ~~\{x,y,z\}=- \{y,x,z\},\\
&&(LY3) ~~ \circlearrowleft_{x,y,z}[[x, y], z]+\circlearrowleft_{x,y,z}\{x, y, z\}=0,\\
&&(LY4) ~~\circlearrowleft_{x,y,z}\{[x, y], z, a\} =0,\\
&&(LY5) ~~\{a, b, [x, y]\}=[\{a, b, x\}, y]+[x,\{a, b, y\}],\\
&&(LY6) ~~ \{a, b, \{x, y, z\}\}=\{\{a, b, x\}, y, z\}+ \{x,  \{a, b, y\}, z\}+ \{x,  y, \{a, b, z\}\},
\end{eqnarray*}
for   $ x, y, z, a, b\in L$, and  where $\circlearrowleft_{x,y,z}$ denotes the summation over the cyclic permutations of $x,y,z$, that is $\circlearrowleft_{x,y,z}[[x,y],z]=[[x,y],z]+[[z,x],y]+[[y,z],x]$.
\end{definition}

\begin{example} \label{exam:Lie algebra}
 Let $(L, [-, -])$  be a Lie algebra.  Define a  ternary bracket on $L$ by
 $$\{x,y,z\}=[[x,y],z],\ \ \  \  \forall x,y,z\in L.$$
 Then  $(L, [-, -], \{-, -, -\})$ becomes  a     Lie-Yamaguti  algebra.
\end{example}

\begin{example}
 Let $(L, [-, -])$  be a Lie algebra  with a reductive decomposition $L=N\oplus M$, i.e.  $[N,M]\subseteq N$ and $[N,M]\subseteq M$.
Define bilinear bracket $[-,-]_M$ and trilinear bracket $\{-,-,-\}_M$ on $M$ by
the projections of the Lie bracket:
$$[x,y]_M=\pi_M([x,y]),\{x,y,z\}_M=[\pi_N([x,y]),z], \forall x,y,z\in M,$$
where $\pi_N:L\rightarrow N,\pi_M:L\rightarrow M$ are the projection map.  Then, $(M, [-,-]_M, \{-, -, -\}_M)$ is a   Lie-Yamaguti  algebra.
\end{example}

\begin{example}  \label{exam:Leibniz algebra}
 Let $(L, \star)$  be a Leibniz algebra.   Define a binary and ternary bracket on $L$ by
 $$[x,y]=x\star y-y\star x, \{x,y,z\}=-(x\star y)\star z, \forall x,y,z\in L.$$
 Then $(L, [-, -], \{-, -, -\})$ becomes  a     Lie-Yamaguti  algebra.
\end{example}

\begin{example} \label{exam:2-dimensional LY algebra}
Let $L$ be a two-dimensional  vector space  with a basis $\varepsilon_1,\varepsilon_2$.  If we define  a binary non-zero bracket $[-,-]$ and a ternary non-zero bracket $\{-, -, -\}$ on $L$ as follows:
$$[\varepsilon_1,\varepsilon_2]=-[\varepsilon_2,\varepsilon_1]=\varepsilon_1, \{\varepsilon_1,\varepsilon_2, \varepsilon_2\}=-\{\varepsilon_2,\varepsilon_1, \varepsilon_2\}=\varepsilon_1,$$
then $(L, [-, -], \{-, -, -\})$ is a   Lie-Yamaguti  algebra.
\end{example}

\begin{example} \label{exam:infinite-dimensional LY algebra}
Let $L=C^{\infty}([0,1])$ endowed with the following brackets
\begin{align*}
[f_1,f_2]&=f_1\circ f'_2-f_2\circ f'_1,\\
\{f_1,f_2,f_3\}&=f_1\circ f'_2\circ f'_3-f_2\circ f'_1\circ f'_3-f_3\circ( f_1\circ f''_2-f_2\circ f''_1),
\end{align*}
for  $f_1,f_2,f_3\in L$. Then $(L, [-, -], \{-, -, -\})$ is a Lie-Yamaguti algebra.
\end{example}

%\begin{example} \label{exam:3-dimensional LY algebra}
%Let $\mathfrak{L}$ be a 3-dimensional  vector space  with a basis $\varepsilon_1$, $\varepsilon_2, \varepsilon_3$.  If we define  a binary non-zero bracket $[-,-]$ and a ternary non-zero bracket $\{-, -, -\}$ on $\mathfrak{L}$ %as follows:
%$$[\varepsilon_1,\varepsilon_2]=-[\varepsilon_2,\varepsilon_1]=\varepsilon_3, \{\varepsilon_1,\varepsilon_2, \varepsilon_1\}=-\{\varepsilon_2,\varepsilon_1, \varepsilon_1\}=\varepsilon_3,$$
%then $(\mathfrak{L}, [-, -], \{-, -, -\})$ is a   Lie-Yamaguti  algebra.
%\end{example}

\begin{definition}
Let $(L, [-, -], \{-, -, -\})$  be a  Lie-Yamaguti   algebra   and $V$ be a vector space. A representation of $L$ on $V$ consists of a linear map  $\rho: L\rightarrow \mathfrak{gl}(V)$ and two bilinear maps $\theta:L\times L\rightarrow \mathfrak{gl}(V)$ such that
\begin{align}
&~\theta([x,y],a) =\theta(x,a)\rho(y)-\theta(y,a)\rho(x),\label{2.1}\\
&~D(a,b)\rho(x)=\rho(x)D(a,b)+\rho(\{a,b,x\}),\label{2.2}\\
&~\theta(x,[a,b])=\rho(a)\theta(x,b)-\rho(b)\theta(x,a),\label{2.3}\\
&~D(a,b)\theta(x,y)=\theta(x,y)D(a,b)+\theta(\{a,b,x\},y)+\theta(x,\{a,b,y\}),\label{2.4}\\
&~\theta(a,\{x,y,z\}) =\theta(y,z)\theta(a,x)-\theta(x,z)\theta(a,y)+D(x,y)\theta(a,z),\label{2.5}
\end{align}
for   $x,y,z,a,b\in L$,
where the bilinear map $D:L\times L\rightarrow \mathfrak{gl}(V)$ is given by
\begin{align}
&~D(x,y)=\theta(y,x)-\theta(x,y)-\rho([x,y])+\rho(x)\rho(y)-\rho(y)\rho(x),\label{2.6}
\end{align}
 In this case, we also call $V$ an $L$-module. We denote a representation by $(V; \rho, \theta)$.
\end{definition}
It can be concluded from Eq. \eqref{2.6}  that
\begin{align}
&~D([x,y],z)+D([y,z],x)+D([z,x],y)=0,\label{2.7}\\
&~D(a,b)D(x,y)=D(x,y)D(a,b)+D(\{a,b,x\},y)+D(x,\{a,b,y\}),\label{2.8}
%&~\theta(\{x,y,z\},a)=\theta(x,a)\theta(z,y)-\theta(y,w)\theta(z,x)-\theta(z,a)D(x,y).\label{2.9}
\end{align}

\begin{example}
Let $(L, [-, -], \{-, -, -\})$  be a  Lie-Yamaguti   algebra.  We   define linear maps $\mathrm{ad}: L\rightarrow \mathfrak{gl}(L),\mathcal{L},\mathcal{R}:\otimes^2 L\rightarrow \mathfrak{gl}(L)$ by
\begin{align*}
& \mathrm{ad}(x)(z):=[x,z],\mathcal{R}(x,y)(z):=\{z,x,y\}, ~~ \text{also}, \\
&\mathcal{L}(x,y)(z)=\mathcal{R}(y,x)(z)-\mathcal{R}(x,y)(z)-\mathrm{ad}([x,y])(z)+\mathrm{ad}(x)\mathrm{ad}(y)(z)-\mathrm{ad}(y)\mathrm{ad}(x)(z)\\
&=\{x,y,z\},
\end{align*}for  $x,y,z\in L$.  Then $(L;\mathrm{ad}, \mathcal{R})$ forms a representation of $L$ on itself,  called the adjoint representation.
\end{example}

Let us recall the Yamaguti cohomology theory on Lie-Yamaguti algebras in \cite{Y67}.  Let $(V; \rho, \theta, D)$ be a
representation of a Lie-Yamaguti algebra $L$, and we denote the set of $(n+1)$-cochains by
$\mathcal{C}^{n+1}_{\mathrm{LY}}(L,V)$, where
\begin{equation*}
\mathcal{C}_{\mathrm{LY}}^{n+1}(L,V)= \left\{ \begin{array}{ll}
\mathrm{Hom}(\underbrace{\wedge^2 L\otimes\cdots\otimes\wedge^2 L}_n,V)\times \mathrm{Hom}(\underbrace{\wedge^2 L\otimes\cdots\otimes\wedge^2 L}_n\otimes L,V) &\mbox{ \mbox{}  $n\geq 1,$  }\\
$$\mathrm{Hom}(L,V)$$ &\mbox{ \mbox{}  $ n=0$.  }
 \end{array}  \right.
\end{equation*}

In the sequel, we recall the coboundary map of $(n+1)$-cochains on a Lie-Yamaguti algebra $L$ with the
coefficients in the representation $(V; \rho, \theta, D)$:

If $n\geq 1$, for any $(f,g)\in \mathcal{C}_{\mathrm{LY}}^{n+1}(L,V)$, $ \mathfrak{X}_i=x_i\wedge y_i\in \wedge^2 L, (i=1,2,\cdots,n+1), z\in L$,  the coboundary map $\delta^{n+1}=(\delta^{n+1}_I,\delta^{n+1}_{II}):\mathcal{C}_{\mathrm{LY}}^{n+1}(L,V)\rightarrow \mathcal{C}_{\mathrm{LY}}^{n+2}(L,V), (f,g)\mapsto (\delta^{n+1}_I(f,g),\delta^{n+1}_{II}(f,g))$ is given as follows:
\begin{align*}
&\delta^{n+1}_I(f,g)(\mathfrak{X}_1,\cdots,\mathfrak{X}_{n+1})\\
=&(-1)^n(\rho(x_{n+1})g(\mathfrak{X}_1,\cdots,\mathfrak{X}_{n},y_{n+1})-\rho(y_{n+1})g(\mathfrak{X}_1,\cdots,\mathfrak{X}_{n},x_{n+1})\\
&-g(\mathfrak{X}_1,\cdots,\mathfrak{X}_{n},[x_{n+1},y_{n+1}]))+\sum_{k=1}^{n}(-1)^{k+1}D(\mathfrak{X}_k)f(\mathfrak{X}_1,\cdots,\widehat{\mathfrak{X}_{k}}\cdots,\mathfrak{X}_{n+1})\\
&+\sum_{1\leq k<l\leq n+1}(-1)^k f(\mathfrak{X}_1,\cdots,\widehat{\mathfrak{X}_{k}}\cdots,\{x_k,y_k,x_l\}\wedge y_l+x_l\wedge \{x_k,y_k,y_l\},\cdots,\mathfrak{X}_{n+1}),\\
&\delta^{n+1}_{II}(f,g)(\mathfrak{X}_1,\cdots,\mathfrak{X}_{n+1},z)\\
=&(-1)^n(\theta(y_{n+1},z)g(\mathfrak{X}_1,\cdots,\mathfrak{X}_{n},x_{n+1})-\theta(x_{n+1},z)g(\mathfrak{X}_1,\cdots,\mathfrak{X}_{n},y_{n+1}))\\
&+\sum_{k=1}^{n+1}(-1)^{k+1}D(\mathfrak{X}_k)g(\mathfrak{X}_1,\cdots,\widehat{\mathfrak{X}_{k}}\cdots,\mathfrak{X}_{n+1},z)\\
&+\sum_{1\leq k<l\leq n+1}(-1)^k g(\mathfrak{X}_1,\cdots,\widehat{\mathfrak{X}_{k}}\cdots,\{x_k,y_k,x_l\}\wedge y_l+x_l\wedge \{x_k,y_k,y_l\},\cdots,\mathfrak{X}_{n+1},z)\\
&+\sum_{k=1}^{n+1}(-1)^kg(\mathfrak{X}_1,\cdots,\widehat{\mathfrak{X}_{k}}\cdots,\mathfrak{X}_{n+1},\{x_k,y_k,z\}).
\end{align*}
where $~\widehat{}~$ denotes omission.
For the case that $n=0$, for any $f\in \mathcal{C}_{\mathrm{LY}}^1(L,V)$ , the coboundary map
 $\delta^1=(\delta^1_I,\delta^1_{II})\text{:}$\, $\mathcal{C}_{\mathrm{LY}}^1(L,V)\rightarrow \mathcal{C}_{\mathrm{LY}}^2(L,V),f\rightarrow (\delta^1_I(f),\delta^1_{II}(f))$ is given by:
\begin{align*}
\delta^1_I(f)(x,y)=&\rho(x)f(y)-\rho(y)f(x)-f([x,y]),\\
\delta^1_{II}(f)(x,y,z)=&D(x,y)f(z)+\theta(y,z)f(x)-\theta(x,z)f(y)-f(\{x,y,z\}).
\end{align*}
The corresponding cohomology groups are denoted by $\mathcal{H}^{\bullet}_{\mathrm{LY}}(L,V).$

\section{ Representations of Reynolds Lie-Yamaguti algebras}\label{sec:Representations}
\def\theequation{\arabic{section}. \arabic{equation}}
\setcounter{equation} {0}

In this section, we introduce the   concept and   representations  of  Reynolds Lie-Yamaguti algebras of arbitrary weight. We also provide various examples and new constructions.

\begin{definition}
Let  $(L, [-, -], \{-, -, -\})$ be a Lie-Yamaguti  algebra and $\lambda$ a given element of  $\mathbb{K}$.  A linear map $T:L\rightarrow L$ is called a    Reynolds operator  of weight $\lambda$ if $T$ satisfies
\begin{align}
 [Tx, Ty]=&T([Tx,y]+[x,Ty]+\lambda[Tx,Ty])\label{3.1}\\
\{Tx,Ty,Tz\}=&T(\{x,Ty,Tz\}+\{Tx,y,Tz\}+\{Tx,Ty,z\}+2\lambda\{Tx,Ty,Tz\}),\label{3.2}
\end{align}
for $ x, y, z\in  L$. Moreover, a  Lie-Yamaguti algebra $L$ with a Reynolds operator $T$ of weight $\lambda$  is called a
 Reynolds   Lie-Yamaguti algebra of weight $\lambda$.  We denote it by $(L, [-, -],  \{-, -, -\}, $ $T)$  or  simply by $(L, T)$.
\end{definition}

\begin{definition}
A homomorphism between two  Reynolds   Lie-Yamaguti algebras  $(L, T)$ and $(L', T')$ of same weight $\lambda$ is a linear map $\phi: L\rightarrow L'$ satisfying
\begin{eqnarray*}
&&\phi[x, y]=[\phi(x), \phi(y)]', ~\phi\{x, y, z\}=\{\phi(x), \phi(y),\phi(z)\}' \   \mathrm{and}\    T'\circ \phi=\phi\circ T,
\end{eqnarray*}
 for  $x,y,z\in L$.
\end{definition}

\begin{remark}
  On the one hand, if a Lie-Yamaguti algebra $(L, [-, -], \{-, -, -\})$   reduces
to a Lie triple system $(L,   \{-, -, -\})$, we obtain the notion of a  Reynolds  Lie triple system  $(L,   \{-, -, -\}, T)$  of weight $(2\lambda)$.
 On the other hand, if a Lie-Yamaguti algebra $(L, [-, -], \{-, -,$ $ -\})$   reduces
to a Lie algebra $(L,  [-, -])$, we get the notion of a  Reynolds Lie algebra $(L,   [-, -], T)$  of weight $\lambda$ immediately.
Thus, all the results given in the sequel can be adapted to the  Reynolds Lie algebra   of weight $\lambda$ and Reynolds  Lie triple system    of weight $\lambda$ context.
\end{remark}

\begin{example}
Any Rota-Baxter operator on Lie-Yamaguti algebras is a  Reynolds operator of weight $\lambda = 0$.
\end{example}

\begin{example}
For any Lie-Yamaguti algebra $L$, the pair $(L, \mathrm{Id}_L)$ is a Reynolds Lie-Yamaguti algebras of  weight -1.
\end{example}

\begin{example}
 Let $(L, T)$  be a   Reynolds  Lie-Yamaguti algebra of  weight $\lambda$.  For   $\lambda'\neq 0\in  \mathbb{K}$, the pair $(L,  \lambda' T)$ is a  Reynolds Lie-Yamaguti algebra of  weight $\frac{\lambda}{\lambda'}$.
\end{example}

\begin{example}
Let $T:L\rightarrow L$ be a    Reynolds operator  of weight $\lambda$ over a Lie algebra $(L, [-, -])$, by $\mathrm{Example}$ \ref{exam:Lie algebra},
we obtain that   $T$ is a    Reynolds operator  of weight $\lambda$ over a Lie-Yamaguti algebra $(L, [-, -], \{-, -, -\})$.
\end{example}

\begin{example}
 Let $(L, \star, T)$  be a Reynolds Leibniz algebra of weight $\lambda$, by $\mathrm{Example}$ \ref{exam:Leibniz algebra}, we obtain that $(L, [-, -], \{-, -, -\}, T)$ is a  Reynolds Lie-Yamaguti algebra of  weight $\lambda$.
\end{example}
\begin{example}
Let $D : L \rightarrow L$ be a derivation on a Lie-Yamaguti algebra $L$. If $(D-\frac{1}{2}\lambda \mathrm{Id}):
L \rightarrow L$ is invertible, then $(L,    (D-\frac{1}{2}\lambda \mathrm{Id})^{-1})$ is a Reynolds  Lie-Yamaguti algebra of  weight $\lambda$.
\end{example}

\begin{example}
Let $(L, [-, -], \{-, -, -\})$ be the two-dimensional Lie-Yamaguti algebra given in $\mathrm{Example}$ \ref{exam:2-dimensional LY algebra}.
Then, for $k_1,k_2,k\neq 0\in \mathbb{K},$
$$T=\left(
        \begin{array}{cc}
          k_1 & k_2  \\
          0 & k
        \end{array}
      \right)$$
is a Reynolds operator of weight $(-\frac{1}{k})$.
\end{example}

\begin{example}
Let $(L, [-, -], \{-, -, -\})$ be the infinite-dimensional Lie-Yamaguti algebra given in $\mathrm{Example}$ \ref{exam:infinite-dimensional LY algebra}.
Then, the integral operator $T:L\rightarrow L$ defined to be
$$T(f)(x)=\int_0^x f(s)ds $$
is a Reynolds operator of weight $\lambda=0$ on $L$.
\end{example}

Next, we introduce the  representation of  Reynolds Lie-Yamaguti algebras of nonzero weight.

\begin{definition}
  A representation of the Reynolds Lie-Yamaguti algebra $(L, T)$  of weight $\lambda$ is a quintuple $(V; \rho,   \theta, T_V)$ such that the following conditions are satisfied:

  (i) $(V; \rho,   \theta)$ is a representation of the  Lie-Yamaguti algebra $L$;

  (ii)  $T_V:V\rightarrow V$  is a linear map satisfying the following equations
\begin{align}
\rho(Tx)T_Vu=&T_V\big(\rho(Tx)u+\rho(x)T_Vu+\lambda\rho(Tx)T_Vu\big),\label{3.3}\\
\theta(Tx, Ty)T_Vu=&T_V\big(\theta(Tx, Ty) u +\theta(Tx, y)T_Vu+\theta(x, Ty)T_Vu+2\lambda\theta(Tx, Ty)T_Vu\big),\label{3.4}
\end{align}
for any $x,y\in \mathfrak{L}$ and $ u\in V.$
\end{definition}

It can be concluded from  Eq. \eqref{3.4}  that
\begin{align}
D(Tx, Ty)T_Vu=&R_V\big(D(Tx, Ty) u +D(Tx, y)T_Vu+D(x, Ty)T_Vu+2\lambda D(Tx, Ty)T_Vu\big),\label{3.5}
\end{align}

\begin{example}
It is easy to  see that  $(L;\mathrm{ad}, \mathcal{R},T)$ is an adjoint representation of the Reynolds Lie-Yamaguti algebra $(L, T)$  of weight $\lambda$.
\end{example}

\begin{example}
Let $(L, T)$  be a Reynolds Lie-Yamaguti algebra   of weight $\lambda$ and  $(V_i; \rho_i,   \theta_i, $ $T_{V_i})_{i\in I}$
be a family of representation of it. Then the quadruple $(\oplus_{i\in I}V_i; (\rho_i)_{i\in I},   (\theta_i)_{i\in I}, \oplus_{i\in I}T_{V_i})$ is a representation of
the Reynolds Lie-Yamaguti algebra $(L, T)$  of weight $\lambda$.
\end{example}

In the following,  we  construct the semidirect product of the  Reynolds Lie-Yamaguti algebra   of weight $\lambda$.

\begin{proposition}
If  $(V; \rho,   \theta,   T_V)$  is a representation of the  Reynolds Lie-Yamaguti algebra $(L, T)$  of weight $\lambda$, then $L \oplus V$ becomes a   Reynolds Lie-Yamaguti algebra    of weight $\lambda$   under the following maps:
\begin{align*}
[x+u, y+v]_{\ltimes}:=&[x, y]+\rho(x)v-  \rho(y)u,\\
\{x+u, y+v, z+w\}_{\ltimes}:=&\{x, y, z\}+D(x, y)w- \theta(x, z)v+ \theta(y, z)u,\\
T\oplus T_V(x+u):=&Tx+T_Vu,
\end{align*}
for all  $x, y, z\in L$ and $u, v, w\in V$. In the case, the Reynolds Lie-Yamaguti algebra  $L \oplus V$  of weight $\lambda$  is called a semidirect product of $L$ and $V$, denoted by $L\ltimes V=(L \oplus V,[-,-]_{\ltimes},\{-,-,-\}_{\ltimes},T\oplus T_V)$.
\end{proposition}

\begin{proof}
For any $x, y, z\in L$ and $u, v, w\in V$,   by  Eqs. \eqref{3.1}--\eqref{3.5},  we have
\begin{align*}
&[T\oplus T_V(x+u), T\oplus T_V(y+v)]_{\ltimes}\\
&\ =[Tx, Ty]+\rho(Tx)(T_Vv)-  \rho(Ty)(T_Vu)\\
&\ =T([Tx,y]+[x,Ty]+\lambda[Tx,Ty])+T_V(\rho(Tx)v+\rho(x)T_Vv+\lambda\rho(Tx)T_Vv)\\
&\ \ \ \ -  T_V(\rho(Ty)u+\rho(y)T_Vu+\lambda\rho(Ty)T_Vu)\\
&\ =(T\oplus T_V)([(T\oplus T_V)(x+u), y+v]_{\ltimes}+[x+u, (T\oplus T_V)(y+v)]_{\ltimes}\\
&\ \ \ \  +\lambda[(T\oplus T_V)(x+u), (T\oplus T_V)(y+v)]_{\ltimes},\\
& [T\oplus T_V(x+u), T\oplus T_V(y+v), T\oplus T_V(z+w)\}_{\ltimes}\\
&\ = \{Tx, Ty, Tz\} +D(Tx, Ty)T_V(w)-\theta(Tx, Tz)T_V(v)+\theta(Ty, Tz)T_V(u)\\
&\ =T\big( \{Tx, Ty, z\}+\{x, Ty, Tz\}+\{Tx, y, Tz\}+2\lambda\{Tx, Ty, Tz\}\big)\\
&\ \ \ \  +T_V\big(D(Tx, Ty)w +D(Tx, y)T_V w +D(x, Ty)T_V w +2\lambda D(Tx, Ty)T_Vw\big)\\
&\ \ \ \ -T_V\big(\theta(Tx, Tz) v +\theta(Tx, z)T_V v +\theta(x, Tz)T_V v +2\lambda\theta(Tx, Tz)T_Vv\big)\\
&\ \ \ \  +T_V\big(\theta(Ty, Tz)u +\theta(Ty, z)T_V(u)+\theta(y, Tz)T_V(u)+2\lambda\theta(Ty, Tz)T_Vu\big)\\
&\ =T\oplus T_V\Big( \{T\oplus T_V(x+u), T\oplus T_V(y+v), z+w\}_{\ltimes}+\{x+u, T\oplus T_V(y+v),\\
&\ \ \ \ \ \ T\oplus T_V(z+w)\}_{\ltimes}+\{T\oplus T_V(x+u), y+v, T\oplus T_V(z+w)\}_{\ltimes}\\
&\ \ \ \  +2\lambda\{T\oplus T_V(x+u), T\oplus T_V(y+v),T\oplus T_V(z+w)\}_{\ltimes}\Big).
\end{align*}
Therefore, $(L \oplus V,[-,-]_{\ltimes},\{-,-,-\}_{\ltimes},T\oplus T_V)$ is a  Reynolds Lie-Yamaguti algebra   of weight $\lambda$.
\end{proof}

\begin{proposition}\label{prop:nLY}
Let $(L,   T)$  be  a    Reynolds  Lie-Yamaguti algebra   of weight $\lambda$.  Define   new operations as follows:
\begin{align}
 [x, y]_T=&[Tx,y]+[x,Ty]+\lambda[Tx, Ty],\label{3.6}\\
\{x,y,z\}_T=&\{x,Ty,Tz\}+\{Tx,y,Tz\}+\{Tx,Ty,z\}+2\lambda\{Tx, Ty, Tz\}, \label{3.7}
\end{align}
for $x,y,z\in L.$ Then,
 \item[(i)] the triple $(L,  [-, -]_T,  \{-, -, -\}_T)$ is a    new Lie-Yamaguti algebra. We denote this Lie-Yamaguti algebra by $L_T.$
 \item[(ii)] the pair $(L_T, T)$ also forms    a  Reynolds  Lie-Yamaguti algebra   of weight $\lambda$.
 \item[(iii)] the map $T:(L,  [-, -]_T,  \{-, -, -\}_T,T)\rightarrow (L,  [-, -],  \{-, -, -\},T )$ is a  morphism of Reynolds  Lie-Yamaguti algebra   of weight $\lambda$.
\end{proposition}

\begin{proof}  (i)  It is directly verified that brackets $[-, -]_T $ and  $ \{-, -, -\}_T$ satisfying Eqs. (LY1)-- (LY6), so the triple $(L,  [-, -]_T,  \{-, -, -\}_T)$ is a   Lie-Yamaguti algebra.
  \item[(ii)]  For any $x,y,z\in L,$ by Eqs. \eqref{3.1}, \eqref{3.2}, \eqref{3.6} and   \eqref{3.7}, we have
  \begin{align*}
& [Tx, Ty]_T\\
& =[T^2x,Ty]+[Tx,T^2y]+\lambda[T^2x, T^2y] \\
& =T\big([Tx,Ty]+[T^2x,y]+\lambda[T^2x, Ty]\big)+T\big([Tx,Ty]+[x,T^2y]+\lambda[Tx, T^2y]\big)\\
&\  \  \ \ +\lambda T\big([Tx,T^2y]+[T^2x,Ty]+\lambda[T^2x, T^2y]\big) \\
& =T\big([x,Ty]_T+[Tx,y]_T+\lambda[Tx, Ty]_T\big), \\
&\{Tx,Ty,Tz\}_T\\
&=\{Tx,T^2y,T^2z\}+\{T^2x,Ty,T^2z\}+\{T^2x,T^2y,Tz\}+2\lambda\{T^2x, T^2y, T^2z\}\\
&=T\big(\{x,T^2y,T^2z\}+\{Tx,Ty,T^2z\}+\{Tx,T^2y,Tz\}+2\lambda\{Tx,T^2y,T^2z\}\big)\\
&\  \  \ \ +T\big(\{Tx,Ty,T^2z\}+\{T^2x,y,T^2z\}+\{T^2x,Ty,Tz\}+2\lambda\{T^2x,Ty,T^2z\}\big)\\
&\  \  \ \ +T\big(\{Tx,T^2y,Tz\}+\{T^2x,Ty,Tz\}+\{T^2x,T^2y,z\}+2\lambda\{T^2x,T^2y,Tz\}\big)\\
&\  \  \ \  +2\lambda T\big(\{Tx, T^2y, T^2z\}+\{T^2x, Ty, T^2z\}+\{T^2x, T^2y, Tz\}+2\lambda\{T^2x, T^2y, T^2z\}\big)\\
&=T\big(\{x,Ty,Tz\}_T+\{Tx,y,Tz\}_T+\{Tx,Ty,z\}_T+2\lambda\{Tx, Ty, Tz\}_T\big),
\end{align*}
which implies that $T$ is a Reynolds operator  of weight $\lambda$ on the Lie-Yamaguti algebra  $L_T.$
\item[(iii)]  From Eqs.  \eqref{3.1} and   \eqref{3.2},  $T$ is a  Lie-Yamaguti algebra homomorphism  from $L_T$ to $L$.
 Moreover, $T$ commutes with itself.
Hence, $T$ is a Reynolds  Lie-Yamaguti algebra   of weight $\lambda$ homomorphism from $(L,  [-, -]_T,  \{-, -, -\}_T,T)$ to $(L,  [-, -],  \{-, -, -\},T)$.
\end{proof}

\begin{proposition}\label{prop:nLYr}
Let $(V; \rho,   \theta,  T_V)$  be a representation of the  Reynolds  Lie-Yamaguti algebra  $(L,   T)$ of weight $\lambda$.
Define   linear maps $\rho_T: L \rightarrow \mathfrak{gl}(V), \theta_T: L\times L \rightarrow \mathfrak{gl}(V)$  by
 \begin{align}
 \rho_T(x)u:=&\rho(Tx)u-T_V(\lambda\rho(Tx)u+\rho(x)u), \label{3.8}\\
\theta_T(x, y)u:=&\theta(Tx, Ty) u-T_V\big(2\lambda\theta(Tx, Ty) u+\theta(Tx,  y) u+\theta(x, Ty)u\big), \label{3.9}
\end{align}
for any $x,y\in L$ and $u\in V.$
It can be concluded from  Eq. \eqref{3.9}  that
 \begin{align}
D_T(x, y)u=&D(Tx, Ty) u-T_V\big(2\lambda D(Tx, Ty) u+D(x, Ty) u+D(Tx, y)u\big). \label{3.10}
\end{align}
Then $(V; \rho_T,   \theta_T)$  is a representation of the     Lie-Yamaguti algebra $L_T$.  Moreover, $(V; \rho_T,   \theta_T,  T_V)$  is a representation of the Reynolds  Lie-Yamaguti algebra   $(L_T,T)$ of weight $\lambda$.
\end{proposition}

\begin{proof}
First, through direct verification, $(V; \rho_T,   \theta_T)$  is a representation of the  Lie-Yamaguti algebra $L_T$.
Further, for any $x, y, z\in L$ and $u, v, w\in V$,   by  Eqs.   \eqref{3.3} and \eqref{3.4},  we have
\begin{align*}
&\rho_T(Tx)T_Vu\\
&\ =\rho(T^2x)T_Vu-T_V\big(\lambda\rho(T^2x)T_Vu+\rho(Tx)T_Vu\big)\\
&\ = T_V\big(\rho(Tx)T_Vu+\rho(T^2x)u+\lambda\rho(T^2x)T_Vu\big)\\
\quad &\ \ \ \  - T_V\Big(\lambda T_V\big(\rho(Tx)T_Vu+\rho(T^2x)u+\lambda\rho(T^2x)T_Vu\big)\\
\quad &\ \ \ \ \ \ \ +T_V\big(\rho(Tx)u+\rho(x)T_Vu+\lambda\rho(Tx)T_Vu\big)\Big)\\
&\ =T_V\big(\rho_T(Tx)u+\rho_T(x)T_Vu+\lambda\rho_T(Tx)T_Vu\big),\\
& \theta_T(Tx, Ty)T_{V}v\\
\quad &= \theta(T^2x, T^2y)T_{V}v-T_V\big(2\lambda\theta(T^2x, T^2y)T_{V}v+\theta(Tx,  T^2y)T_{V}v+\theta(T^2x, Ty)T_{V}v\big)\\
\quad &= T_{V}\big(\theta(T^2x, T^2y) v+\theta(T^2x, Ty)T_{V}v+\theta(Tx, T^2y)T_{V}v+2\lambda\theta(T^2x, T^2y)T_{V}v\big)\\
\quad &\ \ \ \ - T_V\Big(2\lambda T_{V}\big(\theta(T^2x, T^2y)v+\theta(T^2x, Ty)T_{V}v+\theta(Tx, T^2y)T_{V}v+2\lambda\theta(T^2x, T^2y)T_{V}v\big)\\
\quad &\ \ \ \ \  \ +T_{V}\big(\theta(Tx,  T^2y)v+\theta(Tx,  Ty)T_{V}v+\theta(x,  T^2y)T_{V}v+2\lambda\theta(Tx,  T^2y)T_{V}v\big)\\
\quad &\ \ \ \ \  \ +T_{V}\big(\theta(T^2x, Ty)v+\theta(T^2x, y)T_{V}v+\theta(Tx, Ty)T_{V}v+2\lambda\theta(T^2x, Ty)T_{V}v\big)\Big)\\
\quad &= T_V\big(\theta_T(Tx, Ty)v+\theta_T(Tx, y)T_{V}v+\theta_T(x,  Ty)T_{V}v+2\lambda\theta_T(Tx, Ty)T_{V}v\big).
\end{align*}
Hence, $(V; \rho_T,   \theta_T, T_V)$  is a representation of the Reynolds  Lie-Yamaguti algebra   $(L_T,T)$ of weight $\lambda$.
\end{proof}
\begin{example}
$(L;\mathrm{ad}_T, \mathcal{R}_T,T)$ is an adjoint representation of the Reynolds  Lie-Yamaguti algebra   $(L_T, T)$ of weight $\lambda$,
where
\begin{align*}
\mathrm{ad}_T(x)(z):=&[Tx,z]-T(\lambda[Tx,z]+[x,z]),\\
\mathcal{R}_T(x,y)(z):=&\{z,Tx,Ty\}-T(2\lambda\{z,Tx,Ty\}+\{z,Tx,y\}+\{z,x,T y\}), ~~\text{also},\\
\mathcal{L}_T(x,y)(z):=&\{Tx,Ty,z\}-T(2\lambda\{Tx,Ty,z\}+\{Tx, y,z\}+\{x,Ty,z\}),
\end{align*}for  $x,y,z\in L_T$.
\end{example}

\section{ Cohomology of Reynolds  Lie-Yamaguti algebras}\label{sec:Cohomology}
\def\theequation{\arabic{section}. \arabic{equation}}
\setcounter{equation} {0}

In this section, we will construct the cohomology of  Reynolds  Lie-Yamaguti algebras of any weight.
In the next section, we will  use the second cohomology group to study the formal deformation and abelian extension of   Reynolds  Lie-Yamaguti algebras.

 Firstly, the cohomology of the Reynolds    operator $T$ of weight $\lambda$  is given by  Yamaguti cohomology \cite{Y67}.

Let $(V; \rho,   \theta,   T_V)$  be a representation of the  Reynolds  Lie-Yamaguti algebra  $(L, T)$ of weight $\lambda$.
Recall that Proposition \ref{prop:nLY} and Proposition \ref{prop:nLYr} give a new Lie-Yamaguti algebra $L_T$ and
a new representation $(V; \rho_T,   \theta_T)$ over $L_T$. Consider the cochain complex of $L_T$ with coefficients in $(V; \rho_T,   \theta_T)$:
\begin{equation*}
(\mathcal{C}_{\mathrm{LY}}^{\bullet}(L_T,V),\partial^\bullet)=(\oplus_{n=0}^{\infty}\mathcal{C}_{\mathrm{LY}}^{n+1}(L_T,V),\partial^\bullet).
\end{equation*}
More precisely,
\begin{equation*}
\mathcal{C}_{\mathrm{LY}}^{n+1}(L_T,V)= \left\{ \begin{array}{ll}
\mathrm{Hom}(\underbrace{\wedge^2 L_T\otimes\cdots\otimes\wedge^2 L_T}_n,V)\times \mathrm{Hom}(\underbrace{\wedge^2 L_T\otimes\cdots\otimes\wedge^2 L_T}_n\otimes L_T,V) &\mbox{ \mbox{}  $n\geq 1,$  }\\
$$\mathrm{Hom}(L_T,V)$$ &\mbox{ \mbox{}  $ n=0$  }
 \end{array}  \right.
\end{equation*}
 and its coboundary map $\partial^{n+1}=(\partial^{n+1}_I,\partial^{n+1}_{II}):\mathcal{C}_{\mathrm{LY}}^{n+1}(L_T,V)\rightarrow \mathcal{C}_{\mathrm{LY}}^{n+2}(L_T,V), (f,g)\mapsto (\partial^{n+1}_I(f,g),\partial^{n+1}_{II}(f,g))$ is given as follows:
\begin{align*}
&\partial^{n+1}_I(f,g)(\mathcal{X}_1,\cdots,\mathcal{X}_{n+1})\\
=&(-1)^n(\rho_T(x_{n+1})g(\mathcal{X}_1,\cdots,\mathcal{X}_{n},y_{n+1})-\rho_T(y_{n+1})g(\mathcal{X}_1,\cdots,\mathcal{X}_{n},x_{n+1})\\
&-g(\mathcal{X}_1,\cdots,\mathcal{X}_{n},[x_{n+1},y_{n+1}]_T))+\sum_{k=1}^{n}(-1)^{k+1}D_T(\mathcal{X}_k)f(\mathcal{X}_1,\cdots,\widehat{\mathcal{X}_{k}}\cdots,\mathcal{X}_{n+1})\\
&+\sum_{1\leq k<l\leq n+1}(-1)^k f(\mathcal{X}_1,\cdots,\widehat{\mathcal{X}_{k}}\cdots,\{x_k,y_k,x_l\}_T\wedge y_l+x_l\wedge \{x_k,y_k,y_l\}_T,\cdots,\mathcal{X}_{n+1}),
\end{align*}
\begin{align*}
&\partial^{n+1}_{II}(f,g)(\mathcal{X}_1,\cdots,\mathcal{X}_{n+1},z)\\
=&(-1)^n(\theta_T(y_{n+1},z)g(\mathcal{X}_1,\cdots,\mathcal{X}_{n},x_{n+1})-\theta_T(x_{n+1},z)g(\mathcal{X}_1,\cdots,\mathcal{X}_{n},y_{n+1}))\\
&+\sum_{k=1}^{n+1}(-1)^{k+1}D_T(\mathcal{X}_k)g(\mathcal{X}_1,\cdots,\widehat{\mathcal{X}_{k}}\cdots,\mathcal{X}_{n+1},z)\\
&+\sum_{1\leq k<l\leq n+1}(-1)^k g(\mathcal{X}_1,\cdots,\widehat{\mathcal{X}_{k}}\cdots,\{x_k,y_k,x_l\}_T\wedge y_l+x_l\wedge \{x_k,y_k,y_l\}_T,\cdots,\mathcal{X}_{n+1},z)\\
&+\sum_{k=1}^{n+1}(-1)^kg(\mathcal{X}_1,\cdots,\widehat{\mathcal{X}_{k}}\cdots,\mathcal{X}_{n+1},\{x_k,y_k,z\}_T),
\end{align*}
where, $n\geq 1$,  $(f,g)\in \mathcal{C}_{\mathrm{LY}}^{n+1}(L_T,V)$, $ \mathcal{X}_i=x_i\wedge y_i\in \wedge^2 L_T, (i=1,2,\cdots,n+1), z\in L_T$.
For   any $f\in \mathcal{C}_{\mathrm{LY}}^1(L_T,V)$, its coboundary map
 $\partial^1=(\partial^1_I,\partial^1_{II})\text{:}$\, $\mathcal{C}_{\mathrm{LY}}^1(L_T,V)\rightarrow \mathcal{C}_{\mathrm{LY}}^2(L_T,V),f\mapsto (\partial^1_I(f),\partial^1_{II}(f))$ is given by:
\begin{align*}
\partial^1_I(f)(x,y)=&\rho_T(x)f(y)-\rho_T(y)f(x)-f([x,y]_T),\\
\partial^1_{II}(f)(x,y,z)=&D_T(x,y)f(z)+\theta_T(y,z)f(x)-\theta_T(x,z)f(y)-f(\{x,y,z\}_T).
\end{align*}

\begin{definition}
Let $(V; \rho,   \theta,  T_V)$  be a representation of the  Reynolds  Lie-Yamaguti algebra   $(L, T)$ of weight $\lambda$.
Then the cochain complex $(\mathcal{C}_{\mathrm{LY}}^{\bullet}(L_T,V),\partial^{\bullet})$ is called the cochain complex
of Reynolds  operator $T$ of weight $\lambda$ with coefficients in $(V; \rho_T,   \theta_T, T_V)$, denoted by $(\mathcal{C}_{\mathrm{RO}}^{\bullet}(L,V),\partial^{\bullet})$.
The cohomology of $(\mathcal{C}_{\mathrm{RO}}^{\ast}(L,V),\partial^{\bullet})$, denoted by $\mathcal{H}_{\mathrm{RO}}^{\ast}(L,V)$, are called the cohomology of Reynolds  operator $T$ of weight $\lambda$ with coefficients in $(V; \rho_T,   \theta_T,T_V)$.
\end{definition}

In particular, when $(L;\mathrm{ad}_T, \mathcal{R}_T,T)$ is the adjoint  representation of $(L_T, T)$, we denote $(\mathcal{C}_{\mathrm{RO}}^{\bullet}(L,L),\partial^{\bullet})$ by $(\mathcal{C}_{\mathrm{RO}}^{\bullet}(L),\partial^{\bullet})$
and call it the cochain complex of Reynolds  operator $T$ of weight $\lambda$, and denote $\mathcal{H}_{\mathrm{RO}}^{\bullet}(L,L)$ by
$\mathcal{H}_{\mathrm{RO}}^{\bullet}(L)$ and call it the cohomology of  Reynolds  operator $T$ of weight $\lambda$.

Next, we will combine the cohomology of Lie-Yamaguti algebras and the cohomology of
Reynolds  operators  of weight $\lambda$ to construct a cohomology theory for Reynolds   Lie-Yamaguti algebras   of weight $\lambda$.

\begin{definition}
Let $(V; \rho,   \theta, T_V)$  be a representation of the  Reynolds   Lie-Yamaguti algebra $(L, T)$ of weight $\lambda$.
For any  $(f,g)\in \mathcal{C}_{\mathrm{LY}}^{n+1}(L,V)$, $n\geq 1$, we define a linear map $\Phi^{n+1}=(\Phi_I^{n+1},\Phi^{n+1}_{II}):$
$\mathcal{C}^{n+1}_{\mathrm{LY}}(L,V)\rightarrow \mathcal{C}^{n+1}_{\mathrm{RO}}(L,V),  (f,g)\mapsto (\Phi_I^{n+1}(f), \Phi_{II}^{n+1}(g))$  by:
\begin{small}
\begin{align*}
\Phi^{n+1}_{I}(f)=&f\circ(T,\cdots,T)\circ((\mathrm{Id}\wedge \mathrm{Id})^{n})-T_V\Big(\sum_{i=1}^{2n}f\circ(T^{i-1}, \mathrm{Id},T^{2n-i})\circ((\mathrm{Id}\wedge \mathrm{Id})^{n})\\
&+(2n-1)\lambda   f\circ(T,\cdots,T)\circ((\mathrm{Id}\wedge \mathrm{Id})^{n})\Big),\\
\Phi^{n+1}_{II}(g)=&g\circ(T,\cdots,T,T)\circ((\mathrm{Id}\wedge \mathrm{Id})^{n}\wedge \mathrm{Id})-T_V\Big(\sum_{i=1}^{2n+1}g\circ(T^{i-1}, \mathrm{Id},T^{2n+1-i})\circ((\mathrm{Id}\wedge \mathrm{Id})^{n}\wedge \mathrm{Id})\\
&+2n \lambda   g\circ(T,\cdots,T,T)\circ((\mathrm{Id}\wedge \mathrm{Id})^{n}\wedge \mathrm{Id})\Big).
\end{align*}
\end{small}
In particular, when $n=0$, define $\Phi^{1}:\mathcal{C}^{1}_{\mathrm{LY}}(L,V)\rightarrow \mathcal{C}^{1}_{\mathrm{RO}}(L,V)$ by
$$\Phi^1(f)=f\circ T-T_V\circ f.$$
\end{definition}

\begin{lemma}\label{lemma:chain map}
The map $\Phi^{n+1}:\mathcal{C}^{n+1}_{\mathrm{LY}}(L,V)\rightarrow \mathcal{C}^{n+1}_{\mathrm{RO}}(L,V)$ is a cochain map, i.e., the following commutative diagram:
$$\aligned
\xymatrix{
 \mathcal{C}^1_{\mathrm{LY}}(L ,V)\ar[r]^-{\delta^1}\ar[d]^-{\Phi^1}&  \mathcal{C}^2_{\mathrm{LY}}(L,V)\ar@{.}[r]\ar[d]^-{\Phi^2}& \mathcal{C}^{n+1}_{\mathrm{LY}}(L,V)\ar[r]^-{\delta^{n+1}}\ar[d]^-{\Phi^{n+1}}& \mathcal{C}^{n+2}_{\mathrm{LY}}(L,V)\ar[d]^{\Phi^{n+2}}\ar@{.}[r]&\\
 \mathcal{C}^1_{{\mathrm{RO}}}(L,V)\ar[r]^-{\partial^1}& \mathcal{C}^2_{{\mathrm{RO}}}(L,V)\ar@{.}[r]&  \mathcal{C}^{n+1}_{{\mathrm{RO}}}(L,V)\ar[r]^-{\partial^{n+1}}& \mathcal{C}^{n+2}_{{\mathrm{RO}}}(L,V)\ar@{.}[r]&.}
 \endaligned$$
\end{lemma}

 \begin{proof}
For any $(f,g)\in \mathcal{C}_{\mathrm{LY}}^{n+1}(L,V)$, we have
\begin{align}
&\partial_I^{n+1}(\Phi^{n+1}(f,g))(\mathcal{X}_1,\cdots,\mathcal{X}_{n+1})\nonumber\\
=&(-1)^n(\rho_T(x_{n+1})(\Phi^{n+1}_{II}g)(\mathcal{X}_1,\cdots,\mathcal{X}_{n},y_{n+1})-\rho_T(y_{n+1})(\Phi^{n+1}_{II}g)(\mathcal{X}_1,\cdots,\mathcal{X}_{n},x_{n+1})\nonumber\\
&-(\Phi^{n+1}_{II}g)(\mathcal{X}_1,\cdots,\mathcal{X}_{n},[x_{n+1},y_{n+1}]_T))+\sum_{k=1}^{n}(-1)^{k+1}D_T(\mathcal{X}_k)(\Phi^{n+1}_{I}f)(\mathcal{X}_1,\cdots,\widehat{\mathcal{X}_{k}}\cdots,\mathcal{X}_{n+1})\nonumber\\
&+\sum_{1\leq k<l\leq n+1}(-1)^k (\Phi^{n+1}_{I}f)(\mathcal{X}_1,\cdots,\widehat{\mathcal{X}_{k}}\cdots,\{x_k,y_k,x_l\}_T\wedge y_l+x_l\wedge \{x_k,y_k,y_l\}_T,\cdots,\mathcal{X}_{n+1}), \label{4.1}
\end{align}
and
\begin{align}
&\Phi_I^{n+2}(\delta_I^{n+1} (f,g))(\mathcal{X}_1,\cdots,\mathcal{X}_{n+1})\nonumber\\
=&\delta_I^{n+1} (f,g)(T\mathcal{X}_1,\cdots,T\mathcal{X}_{n+1})-T_V\sum_{i=1}^{n+1}\delta_I^{n+1} (f,g)(T\mathcal{X}_1,\cdots,x_{i}\wedge T y_{i}+T x_{i}\wedge y_{i},\cdots,T\mathcal{X}_{n+1})\nonumber\\
&-(2n-1)\lambda T_V\circ \delta_I^{n+1} (f,g)(T\mathcal{X}_1,\cdots,T\mathcal{X}_{n+1}). \label{4.2}
\end{align}
Using Eqs.   \eqref{3.1}--\eqref{3.10} and further expanding  Eqs.   \eqref{4.1} and \eqref{4.2},
 we  have, $\Phi_I^{n+2}\circ\delta_I^{n+1}=\partial_I^{n+1}\circ\Phi^{n+1}$. Similarly, we also have  $\Phi_{II}^{n+2}\circ\delta_{II}^{n+1}=\partial_{II}^{n+1}\circ\Phi^{n+1}$.
\end{proof}

\begin{definition}
Let $(V; \rho,   \theta,T_V)$  be a representation of the  Reynolds   Lie-Yamaguti algebra   $(L, T)$ of weight $\lambda$.
We define the cochain complex $(\mathcal{C}_{\mathrm{RLY}}^{\bullet}(L,V),\mathrm{d}^{\bullet})$  of  Reynolds   Lie-Yamaguti algebra   $(L, T)$ of weight $\lambda$ with coefficients in $(V; \rho,   \theta,T_V)$ to the negative shift of the mapping cone of $\Phi^{\bullet}$, that is, let
$$\mathcal{C}^{1}_{\mathrm{RLY}}(L,V)=\mathcal{C}^{1}_{\mathrm{LY}}(L,V) ~~\text{and}~~
\mathcal{C}_{\mathrm{RLY}}^{n+1}(L,V):=
\mathcal{C}^{n+1}_{\mathrm{LY}}(L,V)\oplus \mathcal{C}^{n}_{\mathrm{RO}}(L,V), ~~\text{for}~~n\geq 1,$$
and the  coboundary map  $\mathrm{d}^{n+1}:\mathcal{C}_{\mathrm{RLY}}^{n+1}(L,V)\rightarrow \mathcal{C}_{\mathrm{RLY}}^{n+2}(L,V)$  by
\begin{align*}
&\mathrm{d}^1(f_1)=(\delta^1 f_1, -\Phi^1 (f_1)), ~~\text{for}~~ f_1\in \mathcal{C}_{\mathrm{RLY}}^{1}(L,V);\\
&\mathrm{d}^{2}((f_1,g_1),f_2)=(\delta^{2}(f_1,g_1), -\partial^{1} (f_2)- \Phi^{2}(f_1,g_1)),~~\text{for}~~((f_1,g_1),f_2)\in \mathcal{C}_{\mathrm{RLY}}^{2}(L,V);\\
&\mathrm{d}^{n+1}((f_1,g_1),(f_2,g_2))=(\delta^{n+1}(f_1,g_1), -\partial^{n} (f_2,g_2)- \Phi^{n+1}(f_1,g_1)),
\end{align*}
 for $((f_1,g_1),(f_2,g_2))\in \mathcal{C}_{\mathrm{RLY}}^{n+1}(L,V)$, $n\geq 2$.
\end{definition}

The cohomology of  $(\mathcal{C}_{\mathrm{RLY}}^{\bullet}(L,V),\mathrm{d}^{\bullet})$, denoted by
 $\mathcal{H}_{\mathrm{RLY}}^{\bullet}(L,V)$,
 is called the cohomology of the  Reynolds   Lie-Yamaguti algebra   $(L, T)$ of weight $\lambda$  with coefficients in $(V; \rho,   \theta, T_V)$.

 In particular,  when $(V; \rho,   \theta,  T_V)=(L;\mathrm{ad}, \mathcal{R},T)$, we just denote $(\mathcal{C}_{\mathrm{RLY}}^{\bullet}(L,L),\mathrm{d}^\bullet)$, $\mathcal{H}_{\mathrm{RLY}}^{\bullet}(L,L)$  by
 $(\mathcal{C}_{\mathrm{RLY}}^{\bullet}(L),\mathrm{d}^\bullet)$, $\mathcal{H}_{\mathrm{RLY}}^{\bullet}(L)$ respectively, and call them the cochain complex, the cohomology of Reynolds   Lie-Yamaguti algebra   $(L, T)$ of weight $\lambda$   respectively.

It is obvious that there is a  short exact sequence of cochain complexes:
\begin{align*}
0\rightarrow \mathcal{C}_{\mathrm{RO}}^{\bullet-1}(L,V)\stackrel{}{\longrightarrow}\mathcal{C}_{\mathrm{RLY}}^{\bullet}(L,V)\stackrel{}{\longrightarrow}\mathcal{C}_{\mathrm{LY}}^{\bullet}(L,V)\rightarrow 0.
\end{align*}
It induces a long exact sequence of cohomology groups:
\begin{align*}
\cdots\rightarrow \mathcal{H}_{\mathrm{RLY}}^{p}(L,V)\rightarrow \mathcal{H}_{\mathrm{LY}}^{p}(L,V)\rightarrow \mathcal{H}_{\mathrm{RO}}^{p}(L,V)\rightarrow \mathcal{H}_{\mathrm{RLY}}^{p+1}(L,V)\rightarrow \mathcal{H}_{\mathrm{LY}}^{p+1}(L,V)\rightarrow \cdots.
\end{align*}

 \section{  Formal deformations of  Reynolds   Lie-Yamaguti algebras}\label{sec:deformations}
\def\theequation{\arabic{section}. \arabic{equation}}
\setcounter{equation} {0}

In this section, we consider  the formal deformation  of  Reynolds   Lie-Yamaguti algebra    $(L, T)$ of weight $\lambda$.
Let $\mathbb{K}[[t]]$ be a ring of power series of one variable $t$, and let $L[[t]]$ be the set of formal power series over $L$.  If
$(L, [-, -], \{-, -, -\})$ is a Lie-Yamaguti algebra, then there is a Lie-Yamaguti algebra structure over the
ring $\mathbb{K}[[t]]$ on $L[[t]]$ given by
\begin{align*}
&[\sum_{i=0}^{\infty}x_it^i,\sum_{i=0}^{\infty}y_jt^j]=\sum_{s=0}^{\infty}\sum_{i+j=s}[x_i,y_j]t^s,\\
& \{\sum_{i=0}^{\infty}x_it^i,\sum_{i=0}^{\infty}y_jt^j,\sum_{k=0}^{\infty}z_kt^k\}=\sum_{s=0}^{\infty}\sum_{i+j+k=s}\{x_i,y_j,z_k\}t^s.
\end{align*}

\begin{definition}
A   formal deformation of the  Reynolds   Lie-Yamaguti algebra    $(L, T)$ of weight $\lambda$ is a triple  $ (F_t, G_t, T_t)$  of the forms
$$F_t=\sum_{i=0}^{\infty}F_it^i,~~ G_t=\sum_{i=1}^{\infty}G_it^i,~~T_t= \sum_{i=0}^{\infty}T_it^i,$$
such that the following conditions are satisfied:

(i) $((F_i,G_i),T_i)\in \mathcal{C}^{2}_{\mathrm{RLY}}(L);$

(ii) $F_0=[-, -],G_0=\{-, -, -\}$ and $T_0=T;$

(iii) and $(L[[t]], F_t, G_t, T_t)$  is  a   Reynolds   Lie-Yamaguti algebra    of weight $\lambda$ over $\mathbb{K}[[t]]$.
\end{definition}

Let $ (F_t, G_t, T_t)$  be a formal deformation as above. Then, for any $ x, y, z, a, b\in L$, the following equations must hold:
\begin{align*}
&~F_t(x, y)+F_t(y,  x)=0,~~G_t(x,y,z)+G_t(y,x,z)=0,\\
&~ \circlearrowleft_{x,y,z}F_t(F_t(x, y), z)+\circlearrowleft_{x,y,z}G_t(x,y,z)=0,\\
&~\circlearrowleft_{x,y,z}G_t(F_t(x, y), z, a)=0,\\
&~G_t(a, b, F_t(x, y))=F_t(G_t(a, b, x), y)+F_t(x,G_t(a, b, y)),\\
& ~ G_t(a, b, G_t(x, y, z))=G_t(G_t(a, b, x), y, z)+ G_t(x,  G_t(a, b, y), z)+ G_t(x,  y, G_t(a, b, z)),\\
&~ F_t(T_tx, T_ty)=T_t\big(F_t(T_tx,y)+F_t(x,T_ty)+\lambda F_t(T_tx, T_ty)\big), \\
&~G_t(T_tx,T_ty,T_tz)=T_t\big(G_t(x,T_ty,T_tz)+G_t(T_tx,y,T_tz)+G_t(T_tx,TR_ty,z)+2\lambda G_t(T_tx,T_ty,T_tz)\big).
\end{align*}
Collecting the coefficients of $t^n$, we get that the above equations are equivalent to the following equations.
\begin{align}
&~F_n(x, y)+F_n(y,  x)=0,~~G_n(x,y,z)+G_n(y,x,z)=0,\label{5.1}\\
& \sum_{i=0}^n\circlearrowleft_{x,y,z}F_i(F_{n-i}(x, y), z)+\circlearrowleft_{x,y,z}G_n(x,y,z)=0,\label{5.2}\\
&\sum_{i=0}^n\circlearrowleft_{x,y,z}G_i(F_{n-i}(x, y), z, a)=0,\label{5.3}\\
&\sum_{i=0}^nG_i(a, b, F_{n-i}(x, y))=\sum_{i=0}^nF_i(G_{n-i}(a, b, x), y)+\sum_{i=0}^nF_i(x,G_{n-i}(a, b, y)),\label{5.4}
\end{align}
\begin{align}
& \sum_{i=0}^n G_i(a, b, G_{n-i}(x, y, z))\nonumber\\
&=\sum_{i=0}^n\big(G_i(G_{n-i}(a, b, x), y, z)+ G_i(x,  G_{n-i}(a, b, y), z)+ G_i(x,  y, G_{n-i}(a, b, z))\big),\label{5.5}\\
&~ \sum_{i+j+k=n}F_i(T_jx, T_ky)=\sum_{i+j+k=n}T_i\big(F_j(T_kx,y)+F_j(x,T_ky)\big)+\lambda \sum_{i+j+k+l=n}T_i(F_j(T_kx, T_ly)), \label{5.6}\\
&\sum_{i+j+k+l=n}G_i(T_jx,T_ky,T_lz)=\sum_{i+j+k+l=n}T_i\big(G_j(x,T_ky,T_lz)+G_j(T_kx,y,T_lz)+G_j(T_kx,T_ly,z)\big)\nonumber\\
&~~~~~~~~~~~ +2\lambda\sum_{i+j+k+l+m=n}T_i\big(G_j(T_kx,T_ly,T_mz)\big).\label{5.7}
\end{align}
Note that for $n=0$, equations \eqref{5.1}--\eqref{5.7}  are equivalent to $(L, F_0, G_0, T_0)$  is   a   Reynolds   Lie-Yamaguti algebra    of weight $\lambda$.

\begin{proposition}\label{prop:2-cocycle}
Let $ (F_t, G_t, T_t)$  be a formal deformation of  a   Reynolds   Lie-Yamaguti algebra    $(L, T)$ of weight $\lambda$.
Then $((F_1,G_1),T_1)$ is a 2-cocycle in the cochain complex $(\mathcal{C}^{\bullet}_{\mathrm{RLY}}(L),\mathrm{d}^\bullet)$.
\end{proposition}
\begin{proof}
For  $n=1$,   equations \eqref{5.1}--\eqref{5.7} become
\begin{align}
&F_1(x, y)+F_1(y,  x)=0,~~G_1(x,y,z)+G_1(y,x,z)=0,\label{5.8}\\
&[F_{1}(x, y), z]+F_1([x, y], z)+[F_{1}(z, x), y]+F_1([z, x], y)+[F_{1}(y, z), x]+F_1([y, z], x)\nonumber\\
&+G_1(x,y,z)+G_1(z,x,y)+G_1(y,z, x)=0,\label{5.9}\\
&G_1([x, y], z, a)+\{F_{1}(x, y), z, a\}+G_1([z, x], y, a)+\{F_{1}(z, x), y, a\}+G_1([y, z], x, a)\nonumber\\
&+\{F_{1}(y, z), x, a\}=0,\label{5.10}\\
&G_1(a, b, [x, y])+\{a, b, F_{1}(x, y)\}\nonumber\\
&=F_1(\{a, b, x\}, y)+[G_{1}(a, b, x), y]+F_1(x,\{a, b, y\})+[x,G_{1}(a, b, y)],\label{5.11}\\
&  G_1(a, b, \{x, y, z\})+  \{a, b, G_{1}(x, y, z)\}=G_1(\{a, b, x\}, y, z)+ \{G_{1}(a, b, x), y, z\}\nonumber\\
&+ G_1(x,  \{a, b, y\}, z)+ \{x,  G_{1}(a, b, y), z\}+ G_1(x,  y, \{a, b, z\})+ \{x,  y, G_{1}(a, b, z)\}.\label{5.12}
\end{align}
From Eqs. \eqref{5.8}--\eqref{5.12}, we get $(F_1,G_1)\in \mathcal{C}^{2}_{\mathrm{LY}}(L)$ and $\delta^{2}(F_1,G_1)=0$.
\begin{small}
\begin{align}
&F_1(Tx, Ty)+[T_1x, Ty]+[Tx, T_1y]=T_1([Tx,y]+[x,Ty])+T(F_1(Tx,y)+F_1(x,Ty))\nonumber\\
&+T([T_1x,y]+[x,T_1y])+\lambda \Big( T_1[Tx, Ty]+T(F_1(Tx, Ty))+T[T_1x, Ty]+T[Tx, T_1y]\Big), \label{5.13}
\end{align}
\begin{align}
&G_1(Tx,Ty,Tz)+\{T_1x,Ty,Tz\}+\{Tx,T_1y,Tz\}+\{Tx,Ty,T_1z\}=T_1\big(\{x,Ty,Tz\}+\{Tx,y,Tz\}\nonumber\\
&+\{Tx,Ty,z\}\big)+T\big(G_1(x,Ty,Tz)+G_1(Tx,y,Tz)+G_1(Tx,Ty,z)\big)+T\big(\{x,T_1y,Tz\}\nonumber\\
&+\{T_1x,y,Tz\}+\{T_1x,Ty,z\}+\{x,Ty,T_1z\}+\{Tx,y,T_1z\}+\{Tx,T_1y,z\}\big)+2\lambda\Big(T_1\{Tx,Ty,Tz\}\nonumber\\
&+TG_1(Tx,Ty,Tz)+T\{T_1x,Ty,Tz\}+T\{Tx,T_1y,Tz\})+T\{Tx,Ty,T_1z\}\Big).\label{5.14}
\end{align}
\end{small}
 Further from    Eqs. \eqref{5.13} and \eqref{5.14}, we have
$-\partial_I^{1} (T_1)- \Phi_I^{2}(F_1)=0$ and $-\partial_{II}^{1} (T_1)- \Phi_{II}^{2}(G_1)=0$  respectively.
 Hence, $\mathrm{d}^{2}((F_1,G_1),T_1)=0,$ that is, $((F_1,G_1),T_1)$ is a 2-cocycle in   $(\mathcal{C}^{\bullet}_{\mathrm{MRBLY}}(L),\mathrm{d}^\bullet)$.
\end{proof}

\begin{definition}
The 2-cocycle  $((F_1,G_1),T_1)$ is called the infinitesimal of the   formal deformation $ (F_t, G_t, T_t)$ of a Reynolds   Lie-Yamaguti algebra    $(L, T)$  of weight $\lambda$.
\end{definition}

\begin{definition}
Let $(F_t, G_t, T_t)$ and $(F'_t, G'_t, T'_t)$ be two formal deformations of
a Reynolds   Lie-Yamaguti algebra    $(L, T)$ of weight $\lambda$. A formal isomorphism from $(L[[t]], F_t, G_t, T_t)$ to $(L[[t]],F'_t, G'_t, T'_t)$ is a power series
$\varphi_t=\sum_{i=o}^{\infty}\varphi_it^i:L[[t]]\rightarrow L[[t]]$ , where  $\varphi_i:L\rightarrow L$ are linear maps
with $\varphi_0=\mathrm{Id}_L$, such that:
\begin{align}
& \varphi_t \circ F_t=F'_t\circ( \varphi_t\otimes \varphi_t),\label{5.15}\\
&  \varphi_t \circ G_t=G'_t\circ( \varphi_t\otimes \varphi_t\otimes \varphi_t),\label{5.16}\\
&  \varphi_t \circ T_t=T'_t\circ  \varphi_t,\label{5.17}
\end{align}
In this case, we say that the two   formal deformations $(F_t, G_t, T_t)$ and $(F'_t, G'_t, T'_t)$
are equivalent.
\end{definition}

\begin{proposition}
The infinitesimals of two equivalent   formal deformations of   $(L,T)$
are in the same cohomology class in $\mathcal{H}^{2}_{\mathrm{RLY}}(L)$.
\end{proposition}
\begin{proof}
Let $\varphi_t:(L[[t]], F_t, G_t, T_t)\rightarrow (L[[t]],F'_t, G'_t, T'_t)$ be a formal isomorphism.
 By expanding Eqs. \eqref{5.15}--\eqref{5.17} and comparing the coefficients of $t$ on both sides, we have
 \begin{align*}
F_1-F'_1=&F_0\circ (\varphi_1\o \mathrm{Id}_L)+F_0\circ (\mathrm{Id}_L \o\varphi_1)-\varphi_1\circ F_0\\
=&\delta_I^1(\varphi_1),\\
G_1-G'_1=&G_0\circ (\varphi_1\o \mathrm{Id}_L\o \mathrm{Id}_L)+G_0\circ (\mathrm{Id}_L\o  \mathrm{Id}_L \o\varphi_1)+G_0\circ (\mathrm{Id}_L\o \varphi_1 \o \mathrm{Id}_L)-\varphi_1\circ G_0\\
=&\delta_{II}^1(\varphi_1),\\
T_{1}-T'_1=& T \circ \varphi_1- \varphi_{1}\circ T\\
=&-\Phi^1(\varphi_1),
 \end{align*}
that is, we have
  $$((F_1,G_1), T_1)-((F'_1,G'_1), T'_1)=(\delta^1(\varphi_1),-\Phi^1(\varphi_1))=\mathrm{d}^1(\varphi_1)\in \mathcal{B}_{\mathrm{RLY}}^{2}(L). $$
  Therefore,  $((F'_1,G'_1), T'_1)$ and $((F_1,G_1), T_1)$ are in the same cohomology class in $\mathcal{H}^{2}_{\mathrm{RLY}}(L)$.
\end{proof}

\begin{definition}
A  formal deformation $(F_t, G_t, T_t)$  of  Reynolds   Lie-Yamaguti algebra $(L,T)$  of weight $\lambda$ is said to be trivial if the deformation $(F_t, G_t, T_t)$ is equivalent
to the undeformed one $(F_0, G_0, T)$.
\end{definition}

\begin{definition}
A  Reynolds   Lie-Yamaguti algebra $(L,T)$  of weight $\lambda$  is said to be rigid if every   formal deformation of $L$ is trivial deformation.
\end{definition}

\begin{theorem}
 Let $(L,T)$  be a  Reynolds   Lie-Yamaguti algebra  of weight $\lambda$ with $\mathcal{H}^{2}_{\mathrm{RLY}}(L)=0$, then
$(L,T)$ is rigid.
\end{theorem}
\begin{proof}
Let $(F_t, G_t, T_t)$ be a   formal deformation of $(L,T)$. From Proposition \ref{prop:2-cocycle},
$((F_1, G_1), T_1)$ is a 2-cocycle. By $\mathcal{H}^{2}_{\mathrm{RLY}}(L)=0$, there exists a 1-cochain
$$\varphi_1 \in \mathcal{C}^1_{\mathrm{RLY}}(L)=\mathcal{C}^1_{\mathrm{LY}}(L)$$
such that $((F_1,G_1), T_1)=  \mathrm{d}^1(\varphi_1), $  that is, $F_1=\delta_I^1(\varphi_1)$, $G_1=\delta_{II}^1(\varphi_1)$ and $R_1=-\Phi^1(\varphi_1)$.

Setting $\varphi_t = \mathrm{id}_L -\varphi_1t$, we get a deformation  $(F'_t, G'_t, T'_t)$, where
 \begin{align*}
F'_t=&\varphi_t^{-1}\circ F_t\circ (\varphi_t\otimes \varphi_t),\\
G'_t=&\varphi_t^{-1}\circ G_t\circ (\varphi_t\otimes \varphi_t\otimes \varphi_t),\\
T'_t=&\varphi_t^{-1}\circ T_t\circ \varphi_t.
 \end{align*}
 It can be easily verify that $F'_1=0,G'_1=0,  T'_1=0$. Then
    $$\begin{array}{rcl} F'_t&=& F_0+F'_2t^2+\cdots,\\
    G'_t&=& G_0+G'_2t^2+\cdots,\\
 T'_t&=& T+T'_2t^2+\cdots.\end{array}$$
  By Eqs.     \eqref{5.1}--\eqref{5.7},
  we see that $((F'_2, G'_2),  T'_2)$  is still a 2-cocyle, so by induction, we can show
that $(F_t, G_t, T_t)$  is equivalent to the trivial deformation $(F_0, G_0, T)$.  Therefore,  $(L,T)$ is rigid.
\end{proof}

  \section{   Abelian extensions of Reynolds   Lie-Yamaguti algebras}\label{sec:extensions}
\def\theequation{\arabic{section}. \arabic{equation}}
\setcounter{equation} {0}

Motivated by the   extensions of Lie-Yamaguti algebras \cite{Z15}, in this section, we consider  abelian extensions of   Reynolds   Lie-Yamaguti algebras  of any weight
and show that they are classified by the second cohomology, as one would expect of a good cohomology theory.

Notice that a vector space $V$ together with a linear map $T_V:V\rightarrow V$ is naturally an
abelian Reynolds   Lie-Yamaguti algebra   of weight $\lambda$   where the multiplications on  $V$ is defined to be $[-,-]_V=0,\{-,-,-\}_V=0.$

\begin{definition}
An abelian extension of  $(L,T)$ by  $(V, T_V)$
 is  a short exact sequence of   Reynolds   Lie-Yamaguti algebras   of weight $\lambda$
$$\begin{CD}%\label{dia:ext}
0@>>> {(V, T_V)} @>i >> (\hat{L},[-, -]_{\hat{L}}, \{-, -, -\}_{\hat{L}},\hat{T}) @>p >> (L,T) @>>>0,
\end{CD}$$
that is, there exists a commutative diagram:
\begin{align*}
\begin{CD}%\label{dia:ext}
0@>>> {V} @>i >> \hat{L} @>p >> L @>>>0\\
@. @V {T_V} VV @V \hat{T} VV @V T VV @.\\
0@>>> {V} @>i >> \hat{L} @>p >> L @>>>0,
\end{CD}
\end{align*}
where    $V$ is an abelian ideal of $\hat{L},$  i.e.,  $[u, v]_{\hat{L}}=0,\{-, u,v\}_{\hat{L}}=\{u,v,-\}_{\hat{L}}=0, \forall u,v\in V$.

We will call $(\hat{L}, \hat{T})=(\hat{L},[-, -]_{\hat{L}}, \{-, -, -\}_{\hat{L}},\hat{T})$ an abelian extension of $(L,T)$ by  $(V, T_V).$
\end{definition}

\begin{definition}
 A   section  of an abelian extension $(\hat{L}, \hat{T})$ of $(L,T)$  by  $(V, T_V)$ is a linear map $s:L\rightarrow \hat{L}$ such that   $p\circ s=\mathrm{id}_L$.
\end{definition}

\begin{definition}
   Let $(\hat{L}_1, \hat{T}_1)$ and  $(\hat{L}_2, \hat{T}_2)$  be two abelian extensions of $(L, T)$  by  $(V, T_V)$. They are said to be  equivalent if  there is an  isomorphism of  Reynolds   Lie-Yamaguti algebras    $\phi:(\hat{L}_1, \hat{T}_1)\rightarrow (\hat{L}_2, \hat{T}_2)$ of weight $\lambda$
such that the following diagram is  commutative:
\begin{align}
\begin{CD}%\label{dia:ext}
0@>>> {(V, T_V)} @>i_1 >> (\hat{L}_1,   \hat{T}_1) @>p_1 >> (L,T) @>>>0\\
@. @| @V \phi VV @| @.\\
0@>>> {(V, T_V)} @>i_2 >> (\hat{L}_2,   \hat{T}_2) @>p_2 >> (L,T) @>>>0.\label{6.1}
\end{CD}
\end{align}
\end{definition}

Now for an  abelian extension $(\hat{L}, \hat{T})$ of $(L,T)$  by  $(V, T_V)$ with a section $s:L\rightarrow\hat{L}$, we define  linear maps $\rho: L \rightarrow \mathfrak{gl}(V)$  and  $\theta: L\times L \rightarrow \mathfrak{gl}(V)$  by
\begin{align*}
&\rho(x)u:=[s(x),i(u)]_{\hat{L}},\\
&\theta(x,y)u:=\{i(u),s(x),s(y)\}_{\hat{L}},  \quad \forall x,y\in L, u\in V.
\end{align*}
In particular, $D(x,y)u=\theta(y,x)u-\theta(x,y)u-\rho([x,y])u+\rho(x)\rho(y)u-\rho(y)\rho(x)u
=\{s(x),s(y),i(u)\}_{\hat{L}}.$
\begin{proposition} \label{prop:representation}
  With the above notations, $(V; \rho, \theta, T_V)$ is a representation of the  Reynolds   Lie-Yamaguti algebra $(L,T)$  of weight $\lambda$.
\end{proposition}
\begin{proof}
For any $x,y\in L$ and $ u\in V,$ ${\hat{T}}s(x)-s(Tx)\in V$ means that $\rho({\hat{T}}s(x))u=\rho(s(Tx))u,\theta({\hat{T}}s(x), {\hat{T}}s(y)) u=\theta(s(Tx), s(Ty)) u$. Therefore,
we have
 \begin{align*}
&\rho(Tx)T_Vu=[s(Tx),i(T_Vu)]_{\hat{L}}=[\hat{T}s(x),\hat{T}i (u)]_{\hat{L}}\\
&=\hat{T}\big([\hat{T}s(x),i (u)]_{\hat{L}}+[s(x),\hat{T}i (u)]_{\hat{L}}+\lambda[\hat{T}s(x),T_Vi (u)]_{\hat{L}}\big)\\
&=T_V\big(\rho(Tx)u+\rho(x)T_Vu+\lambda\rho(Tx)T_Vu),\\
&\theta(Tx,  Ty)T_Vu\\
\quad  &= \{i(T_V u), s(Tx), s(Ty)\}_{\hat{L}}= \{{\hat{T}}i(u), {\hat{T}}s(x), {\hat{T}}s(y)\}_{\hat{L}}\\
\quad  &={\hat{T}}( \{{\hat{T}}i(u), {\hat{T}}s(x), s(y)\}_{\hat{L}}+\{i(u), {\hat{T}}s(x), {\hat{T}}s(y)\}_{\hat{L}}+\{{\hat{T}}i(u), s(x), {\hat{T}}s(y)\}_{\hat{L}}\\
\quad &~~~~ +2\lambda\{{\hat{T}}i(u), {\hat{T}}s(x), {\hat{T}}s(y)\}_{\hat{L}})\\
\quad  &=T_V( \{i(T_V u), s(Tx), s(y)\}_{\hat{L}}+\{i(u), s(Tx), s(Ty)\}_{\hat{L}}+\{i(T_V u), s(x), s(Ty)\}_{\hat{L}}\\
\quad &~~~~ +2\lambda \{i(T_V u), s(Tx), s(Ty)\}_{\hat{L}})\\
\quad  &=T_V\big(\theta(Tx, Ty)u +\theta(Tx, y)T_Vu+\theta(x, Ty)T_Vu+2\lambda\theta(Tx, Ty)T_Vu\big).
\end{align*}
Hence,  $(V; \rho, \theta,T_V)$ is a representation of   $(L,T)$.
\end{proof}

We further define linear maps $\nu:L\times L\rightarrow V$, $\psi:L\times L\times L\rightarrow V$ and $\chi:L\rightarrow V$ respectively by
\begin{align*}
\nu(x,y)&=[s(x), s(y)]_{\hat{L}}-s[x, y],\\
\psi(x,y,z)&=\{s(x), s(y), s(z)\}_{\hat{L}}-s\{x, y,z\},\\
\chi(a)&=\hat{T}s(x)-s(Tx),\quad\forall x,y,z\in L.
\end{align*}
We transfer the the  Reynolds   Lie-Yamaguti algebra  of weight $\lambda$ structure on $\hat{L}$ to $L\oplus V$ by endowing $L\oplus V$ with   multiplications $[-, -]_\nu,\{-,-,-\}_\psi$
and a Reynolds
operator  $T_\chi$  of weight $\lambda$  defined by
\begin{align}
[x+u, y+v]_\nu&=[x, y]+\rho(x)v-\rho(y)u+\nu(x,y), \label{6.2}\\
 \{x+u, y+v, z+w\}_\psi&=\{x, y, z\}+\theta(y, z)u-\theta(x, z)v+D(x, y)w+\psi(x, y, z),\label{6.3}\\
 T_\chi(x+u)&=T(x)+\chi(x)+T_V(u),  \forall x,y,z\in L,\,u,v,w\in V.\label{6.4}
\end{align}

\begin{proposition}\label{prop:2-cocycles}
The 4-tuple $(L\oplus V,[-,-]_\nu, \{-,-,-\}_\psi, T_\chi)$ is a   Reynolds   Lie-Yamaguti algebra  of weight $\lambda$  if and only if
$((\nu,\psi),\chi)$ is a 2-cocycle  of the  Reynolds   Lie-Yamaguti algebra  $(L,T)$  of weight $\lambda$  with the coefficient  in $(V; \rho, \theta, T_V)$.
 In this case,
$$ \begin{CD}%\label{dia:ext}
0@>>> {(V, T_V)} @>i >> (L\oplus V,[-,-]_\nu, \{-,-,-\}_\psi, T_\chi) @>p >>(L,T) @>>>0
\end{CD}$$
 is an abelian extension.
\end{proposition}
\begin{proof}
The map $T_\chi$ is a Reynolds operator    of weight $\lambda$  on $(L\oplus V,[-,-]_\nu, \{-,-,-\}_\psi)$ if and only~if
\begin{align*}
&[T_\chi(x + u), T_\chi(y + v)]_\nu \\
\quad &=T_\chi\big([T_\chi(x + u), y + v]_\nu + [x + u, T_\chi(y + v)]_\nu+\lambda[T_\chi(x + u), T_\chi(y + v)]_\nu\big) , \\
&\{T_\chi(x + u), T_\chi(y + v), T_\chi(z + w)\}_\psi \\
\quad &=T_\chi\big(\{T_\chi(x + u), T_\chi(y + v), z + w\}_\psi + \{x + u, T_\chi(y + v), T_\chi(z + w)\}_\psi \\
\quad & +\{T_\chi(x + u), y + v,T_\chi(z + w)\}_\psi+2\lambda\{T_\chi(x + u), T_\chi(y + v), T_\chi(z + w)\}_\psi\big),
\end{align*}
for   $x,y,z\in L$ and $u,v,w\in V$. Moreover, by   \eqref{3.1}--\eqref{3.10},
the above equations are equivalent to the following equations:
\begin{align}
&\nu(Tx, Ty)+\rho(Tx)\chi(y)-\rho(Ty)\chi(x)=\chi([x,y]_T)+T_V(\rho(x)\chi(y))-T_V(\rho(y)\chi(x))\nonumber\\
&+T_V\big(\nu(Tx,y)+\nu(x,Ty)+\lambda \nu(Tx,Ty)\big)+\lambda T_V\rho(Tx)\chi(y)- \lambda T_V\rho(Ty)\chi(x),\label{6.5}\\
&\psi(Tx, Ty, Tz) + \theta(Ty, Tz)\chi(x)-\theta(Tx, Tz)\chi(y) + D(Tx, Ty)\chi(z)\nonumber\\
\quad &=\chi(\{x, y, z\}_T)+T_V \big(\psi(Tx, Ty, z) + \psi(Tx, y, Tz) + \psi(x, Ty, Tz) +2\lambda \psi(Tx, Ty, Tz)\big)\nonumber\\
 \quad & + T_V \Big(\theta(Ty, z)\chi(x)+\theta(y, Tz)\chi(x)- \theta(Tx, z)\chi(y)-\theta(x, Tz)\chi(y)+D(Tx, y)\chi(z)\nonumber\\
 \quad & +D(x, Ty)\chi(z)+2\lambda\theta(Ty, Tz)\chi(x)-2\lambda\theta(Tx, Tz)\chi(y)+2\lambda D(Tx, Ty)\chi(z)\Big).\label{6.6}
\end{align}
Using Eqs.  \eqref{6.5} and  \eqref{6.6}, we get $-\partial_I^{1} (\chi)- \Phi_I^{2}(\nu)=0$ and $-\partial_{II}^{1} (\chi)- \Phi_{II}^{2}(\psi)=0$  respectively.
Therefore, $ \mathrm{d}^2((\nu,\psi),\chi)=(\delta^2(\nu,\psi),-\partial^1(\chi)-\Phi^{2}(\nu,\psi))=0,$ that is, $((\nu,\psi),\chi)$ is a  2-cocycle.

Conversely, if $((\nu,\psi),\chi)$ is a  2-cocycle  of the  Reynolds   Lie-Yamaguti algebra $(L,T)$  of weight $\lambda$  with the coefficient  in $(V; \rho, \theta, T_V)$,  then we have $ \mathrm{d}^2((\nu,\psi),\chi)=(\delta^2(\nu,\psi),-\partial^1(\chi)-\Phi^{2}(\nu,\psi))=0,$ in which $\delta^2(\nu,\psi)=0$, Eqs.  \eqref{6.5} and \eqref{6.6} hold.
So $(L\oplus V,[-,-]_\nu, \{-,-,$ $-\}_\psi, T_\chi)$   is a   Reynolds   Lie-Yamaguti algebra of weight $\lambda$.
\end{proof}

\begin{proposition}
Let $(\hat{L}, \hat{T})$ be an abelian extension of $(L,T)$  by  $(V, T_V)$ and $s:L\rightarrow\hat{L}$   a section.
\item[(i)]  Different choices of the section $s$ give the same representation on $(V,T_V)$ introduced in Proposition \ref{prop:representation};
\item[(ii)]  If   $((\nu,\psi),\chi)$ is a 2-cocycle   constructed using the section $s$, then  its cohomology class does not depend on the choice of $s$.
\end{proposition}
\begin{proof}
(i) Let $s_1$ and $s_2$ be two distinct sections of $p$. We define a linear map  $\iota: L\rightarrow V$ by $\iota(x)=s_1(x)-s_2(x)$ for $x\in L$.
From $V$ is an abelian ideal of $\hat{L}$, we have
\begin{align*}
 \rho_1(x)u&=[s_1(x),i(u)]_{\hat{L}}=[s_2(x)+\iota(x),i(u)]_{\hat{L}}=[s_2(x),i(u)]_{\hat{L}}=\rho_2(x)u,\\
 \theta_1(x,y)u&=\{i(u),s_1(x),s_1(y)\}_{\hat{L}}=\{i(u),s_2(x)+\iota(x),s_2(y)+\iota(y)\}_{\hat{L}}=\{i(u),s_2(x),s_2(y)\}_{\hat{L}}\\
 &= \theta_2(x,y)u.
\end{align*}
So different choices of $s$ give the same representation on $(V,T_V)$;
\item[(ii)]  Let the 2-cocycles  constructed by  $s_1$ and $s_2$ be $((\nu_1,\psi_1),\chi_1)$ and $((\nu_2,\psi_2),\chi_2)$ respectively.
 Then we have
\begin{align*}
\nu_1(x, y) =& [s_1(x), s_1(y)]_{\hat{L}}-s_1[x, y]\\
\quad =& [s_2(x) +\iota(x), s_2(y) +\iota(y)]_{\hat{L}}- s_2[x, y]-\iota[x, y]\\
\quad =& [s_2(x), s_2(y)]_{\hat{L}} +\rho(x)\iota(y)-\rho(y)\iota(x)-s_2[x, y]- \iota[x, y]\\
\quad =& \nu_2(x, y) + \delta_{I}^1\iota(x, y),\\
\psi_1(x, y, z)  =& \{s_1(x), s_1(y), s_1(z)\}_{\hat{L}}-s_1\{ x, y, z\}\\
\quad =& \{s_2(x) + \iota(x), s_2(y) + \iota(y), s_2(z) + \iota(z)\}_{\hat{L}}- s_2\{x, y, z\}-\lambda\{x, y, z\})\\
\quad =& \{s_2(x), s_2(y), s_2(z)\}_{\hat{L}} + \theta(y, z)\iota(x)-\theta(x, z)\iota(y)\\
\quad &+D(x, y)\iota(z)-s_2\{x, y, z\}- \iota\{x, y, z\}\\
\quad =& \psi_2(x, y, z) + \delta_{II}^1\iota(x, y, z),
\end{align*}
and
\begin{align*}
\chi_1(x) &= \hat{T}s_1(x)- s_1(Tx)\\
\quad &= \hat{T}(s_2(x) + \iota(x))-s_2(Tx) -\iota(Tx)\\
\quad &= \hat{T}s_2(x)-s_2(Tx) + \hat{T}\iota(x)-\iota(Tx)\\
\quad &= \chi_2(x)-\Phi^1\iota(x).
\end{align*}
That is, $((\nu_1,\psi_1),\chi_1)=((\nu_2,\psi_2),\chi_2)+\mathrm{d}^1(\iota)$.  Thus  $((\nu_1,\psi_1),\chi_1)$ and $((\nu_2,\psi_2),\chi_2)$  correspond to the same cohomology
class in $ \mathcal{H}_{\mathrm{RLY}}^{2}(L,V)$.
\end{proof}

\begin{theorem} \label{theorem:classify abelian extensions}
Abelian extensions of a  Reynolds   Lie-Yamaguti algebra $(L,T)$  of weight $\lambda$  by  $(V, T_V)$ are classified by the second cohomology group $\mathcal{H}_{\mathrm{RLY}}^{2}(L,V)$ of $(L,T)$ with coefficients in the representation $(V; \rho, \theta, T_V)$.
\end{theorem}
\begin{proof}
Suppose $(\hat{L}_1, \hat{T}_1)$ and  $(\hat{L}_2, \hat{T}_2)$
 are equivalent abelian extensions   of $(L,T)$  by  $(V, T_V)$  with the associated isomorphism $\phi:(\hat{L}_1, \hat{T}_1)\rightarrow (\hat{L}_2, \hat{T}_2)$ such that the diagram in~\eqref{6.1} is commutative.
 Let $s_1$ be a section of $(\hat{L}_1, \hat{T}_1)$. As $p_2\circ \phi=p_1$, we  have
$$p_2\circ(\phi\circ s_1)=p_1\circ s_1= \mathrm{Id}_L.$$
Thus, $\phi\circ s_1$ is a section of $(\hat{L}_2, \hat{T}_2)$. Denote $s_2:=\phi\circ s_1$. Since $\phi$ is an isomorphism of  Reynolds   Lie-Yamaguti algebras   of weight $\lambda$  such that $\phi|_V=\mathrm{Id}_V$, we have
\begin{align*}
\nu_2(x,y)&=[s_2(x), s_2(y)]_{\hat{L}_2}-s_2[x,y]\\
&=[\phi(s_1(x)), \phi(s_1(y))]_{\hat{L}_2}-\phi(s_1[x, y])\\
&=\phi\big([s_1(x), s_1(y)]_{\hat{L}_1}-s_1[x, y]\big)\\
&=\phi(\nu_1(x,y))\\
&=\nu_1(x,y),\\
\psi_2(x, y, z)&=\{s_2(x), s_2(y), s_2(z)\}_{\hat{L}_2}-s_2\{x, y, z\}\\
\quad &= \phi(\{s_1(x), s_1(y), s_1(z)\}_{\hat{L}_1}-s_1\{x, y, z\})\\
\quad &= \psi_1(x, y, z),\\
\chi_2(x)&=\hat{T}_2s_2(x)-s_2(Tx)\\
\quad &=\hat{T}_2(\phi\circ s_1(x))-\phi\circ s_1(Tx)\\
\quad &=\phi (\hat{T}_1s_1(x)- s_1(Tx))\\
\quad &= \hat{T}_1s_1(x)- s_1(Tx)\\
\quad &=\chi_1(x).
\end{align*}
So, two equivalent abelian extensions give rise to the same element in $\mathcal{H}_{\mathrm{RLY}}^{2}(L,V)$.

Conversely, given two    2-cocycles $((\nu_1,\psi_1),\chi_1)$ and $((\nu_2,\psi_2),\chi_2)$,
we can construct two abelian extensions $(L\oplus V,[-,-]_{\nu_1}, \{-,-,-\}_{\psi_1}, T_{\chi_1})$ and  $(L\oplus V,[-,-]_{\nu_2}, \{-,-,-\}_{\psi_2}, T_{\chi_2})$ via Proposition \ref{prop:2-cocycles}. If they represent the same cohomology class in $\mathcal{H}_{\mathrm{RLY}}^{2}(L,V)$, then  there is  a linear map $\iota: L\rightarrow  V$ such that
 $$((\nu_1,\psi_1),\chi_1)=((\nu_2,\psi_2),\chi_2)+(\delta^1\iota,-\Phi^1\iota).$$

 Define  $\phi_\iota: L\oplus V\rightarrow  L\oplus V$ by
$\phi_\iota(x+u):=x+\iota(x)+u, ~x\in L, u\in V.$  Then we have $\phi_\iota$ is an isomorphism of these two abelian extensions $(L\oplus V,[-,-]_{\nu_1}, \{-,-,-\}_{\psi_1}, T_{\chi_1})$ and  $(L\oplus V,[-,-]_{\nu_2}, \{-,-,-\}_{\psi_2}, T_{\chi_2})$.
\end{proof}

\begin{center}
 {\bf ACKNOWLEDGEMENT}
 \end{center}

  The paper is supported by the NSF of China (No. 12161013) and   Guizhou Provincial Basic Research Program (Natural Science) (No. ZK[2023]025).

\renewcommand{\refname}{REFERENCES}

\end{document}